\documentclass[leqno,12pt]{article}
\usepackage{amsmath,amssymb,amsthm,amsfonts,amscd}
\usepackage{hyperref}
\usepackage[numeric]{amsrefs}
\input colordvi

\newtheorem{thm}[equation]{Theorem}

\newtheorem{lem}[equation]{Lemma}
\newtheorem{rem}[equation]{Remark}

\DeclareMathOperator{\vol}{vol}

\numberwithin{equation}{section}

\newcommand\g{\gamma}

\newcommand\fg{\mathfrak g}

\newcommand\G{\Gamma}

\newcommand{\N}{{\mathbb{N}}}
\newcommand{\Z}{{\mathbb{Z}}}
\newcommand{\R}{{\mathbb{R}}}
\newcommand{\RP}{{\mathbb{RP}}}
\newcommand{\C}{{\mathbb{C}}}

\newcommand{\Q}{{\mathbb{Q}}}
\renewcommand{\O}{{\mathbb{O}}}

\renewcommand\O{{\mathcal O}}
\renewcommand\o{{\mathfrak o}}
\newcommand\E{{\mathcal E}}

\renewcommand\i{^{-1}}
\renewcommand\({\left(}
\renewcommand\){\right)}
\newcommand{\ttwo}[4]{
\(\begin{smallmatrix}{#1} & {#2}
\\ {#3} & {#4} \end{smallmatrix}\)}

\newcommand{\bx}{\hfill$\square$\vspace{.6cm}}

\newcommand{\gobble}[1]{}
  \newcommand{\rangeref}[2]{%
    \ref{#1}--\afterassignment\gobble\fam 0\ref{#2}%
  }

\begin{document}

\title{Pairings of automorphic distributions\thanks{Partially supported by NSF grants  DMS-0500922 and DMS-0601009, and DARPA grant HR0011-04-1-0031}
}

\author{Stephen D. Miller        \and
       Wilfried Schmid 
}

\date{March 28, 2011}

\maketitle

\begin{abstract}
We present a pairing of automorphic distributions that applies in situations where a Lie group acts with an open
orbit on a product of generalized flag varieties.  The pairing gives meaning to an integral of products of
automorphic distributions on these varieties.  This
 generalizes  classical integral representations
or ``Rankin-Selberg integrals'' of $L$-functions, and
 gives new constructions and analytic continuations of automorphic
$L$-functions.

\end{abstract}

\section{Introduction}\label{sec:intro}

Among its many and deep predictions, the Langlands conjectures attach families of $L$-functions to each automorphic
representation. The Langlands $L$-functions are expected to have holomorphic continuations -- or sometimes
meromorphic continuations with poles only at certain specific points -- and to satisfy functional equations. This
has been established in relatively few cases so far, generally by one of two methods. One, the Langlands-Shahidi
method, deduces the meromorphic continuation and functional equation from the analogous properties of Eisenstein
series -- Eisenstein series induced from the automorphic representation in question. The other is called the method
of integral representations, although it is perhaps more a collection of very clever ideas than a systematic
method. These integral representations express the $L$-function as an integral of matrix coefficients (i.e.,
automorphic forms), or products of matrix coefficients, of the automorphic representation; the functional equation
then follows by applying some appropriate involution, which has the desired effect on the argument of the
$L$-function.

The Langlands-Shahidi method presents the functional equation in exactly the form conjectured by Langlands, with
the Gamma factors predicted by him. Its range of applicability is limited, because the group which acts on the
automorphic representation must arise as Levi component of a maximal parabolic in some ambient group. The method of
integral representations typically does not produce the Gamma factors directly. As a result, this method cannot
always rule out all unwanted poles of the $L$-function. The same is true, for different reasons, of the
Langlands-Shahidi method.

In the papers \cites{voronoi,rankin-selberg} we have introduced a new approach: we work not with matrix
coefficients of an automorphic representation, as in the method of integral representations, but with its
automorphic distribution. The latter, in effect, is the datum of the embedding of the automorphic representation
into $L^2(\G\backslash G)$; it does not involve the choice of particular matrix coefficients. Our approach has the
advantage of making the Gamma factors computable and allows us to rule out all the unwanted poles. There is a cost,
of course~-- distributions are more difficult to work with than $C^\infty$ functions.

In this paper we develop an important analytic tool for our program, the analogue, in the setting of automorphic
distributions, of Rankin-Selberg type integrals. Those are integral representations of $L$-functions, as integrals
of products of automorphic forms -- i.e., multilinear pairings of automorphic forms. The analytic properties of the
$L$-functions are then deduced from the integrals in question. Our main result, theorem~\ref{mainthm}, provides an
alternative by establishing an invariant multilinear pairing of automorphic distributions. In our forthcoming paper
\cite{extsquare} we shall use a particular instance of our pairings to obtain new analytic information, such as the
full analytic continuation, of the exterior square $L$-function for $GL(2n)$. We intend to use
theorem~\ref{mainthm} in the future to obtain new results on other $L$-functions, results that are inaccessible by
the Rankin-Selberg and Langlands-Shahidi methods. Indeed, for every integral representation we are aware of, there
exists an analogous distributional pairing of the type covered by theorem~\ref{mainthm}.

\section{Statement of the main theorem}\label{sec:mainthm}

Our main result involves two reductive Lie groups, $G$ and $H$. We suppose that both are finite covers -- typically
but not always the trivial cover -- of the groups of real points of algebraic reductive subgroups $G_\Q \subset
GL(N_G,\Q)$ and $H_\Q \subset GL(N_H,\Q)$, respectively. By an arithmetic subgroup $\G \subset G$, we mean a
subgroup commensurate with the inverse image, in $G$, of the group of integral points $G_\Z \subset GL(N_G,\Z)$. We
use the same terminology in the case of $H$, of course. We shall suppose that $G$ is realized as a subgroup of $H$,
via an inclusion
\begin{equation}
\label{gandh} G\ \hookrightarrow \ H
\end{equation}
that is defined over $\,\Q\,$;\,\ in other words, the inclusion is compatible with a $\,\Q$-homomor\-phism $G_\Q
\to H_\Q$. It then follows that any arithmetic subgroup of $\G_H \subset H$ intersects $G$ in an arithmetic
subgroup $\G_G \subset G$. To simplify various statements, we suppose that
\begin{equation}
\label{cptctr}
Z_H\,,\ \text{the center of $H$\,,\, \ is compact.}
\end{equation}
In the context of automorphic forms on $H$, that is not a restrictive hypothesis: any automorphic form on which the
center of $H$ acts according to a character is completely determined by its restriction to the derived group
$[H,H]$.

We consider a generalized real flag variety $Y$ for $H$ -- a compact real algebraic variety with a transitive
action of $H$, such that the isotropy subgroups are parabolic. If we let $P$ denote the isotropy subgroup at some
base point $y_0 \in Y$, we can make the identification
\begin{equation}
\label{gps2} Y\ \cong \ H/P \,,\ \ \ \ \ y_0 \cong eP\,.
\end{equation}
We do not require that the parabolic subgroup $P\subset H$ be defined over $\,\Q\,$.\,\ Via the embedding
(\ref{gandh}) $G$ acts on $Y$. In addition, we consider the datum of a connected unipotent subgroup $U \subset H$,
which is normalized by $G$ and defined over $\,\Q\,$,\,\ such that\vspace{-4pt}
\begin{equation}
\begin{aligned}
\label{hypo1}
\text{a)}\ \ &\text{$U$ is the unipotent radical of a $\Q$-parabolic subgroup $P_U \subset H$,}
\\[4pt]
\text{b)}\ \ &\text{$G \subset M_U$ for some Langlands decomposition $P_U  = M_U\!\cdot\! A_U \!\cdot \!U$,}
\\[4pt]
\text{c)}\ \ &\text{$U\!\cdot G$\,\ has an open orbit\,\ $\O$\,\ in\,\ $Y$\,, and}
\\[4pt]
\text{d)}\ \ &\text{$G$ preserves the bi-invariant measure on $U$}.
\end{aligned}\vspace{-4pt}
\end{equation}
Since $G$ is reductive and normalizes $U$ it can be extended to a Levi component of $P_U$. Thus b) is implied by
the simpler, but more restrictive condition\vspace{-4pt}
\begin{equation*}
{\rm b}'\,)\ \ \text{$Z_G$\,,\,\ the center of $G$, is compact.}\qquad\qquad\qquad\qquad\qquad\qquad\vspace{-4pt} \tag{\theequation}
\end{equation*}
We note that the semidirect product $\,U\!\cdot G\,$ is defined over $\,\Q\,$,\,\ since both factors are.

At this point two examples may be helpful. In the first, $G$ is the group $SL^{\pm 1}(n,\R)$ equipped with the
standard $\Q$-structure, $H=G\times G \times G$, which contains $G$ diagonally, and $U=\{e\}$. Thus a) and d) are
vacuously satisfied, and the semisimplicity of $G$ implies b). Let $X_n$ denote the flag variety of $G=SL^{\pm
1}(n,\R)$,
\begin{equation}
\label{flagn}
X_n \ = \ \{\,0=F_0 \subset F_1 \subset \dots \subset F_k \subset \dots \subset F_n = \R^n \,\mid \,\dim F_k = k \,\} \,.
\end{equation}
Then $G$ acts on $X_n$ with isotropy subgroups conjugate to the lower triangular subgroup. Note that $\R\mathbb
P^{n-1}$ is a generalized flag variety for $G$. The product $Y =X_n\times X_n\times \R\mathbb P^{n-1}$ can be
regarded as a generalized flag variety for $H$. It is not difficult to see that $G$, via the diagonal embedding
$G\hookrightarrow G\times G$, has a unique open orbit in $X_n\times X_n$, and that the isotropy subgroup of $G$ at
any point of the open orbit is $G$-conjugate to the diagonal Cartan subgroup. The diagonal Cartan has a unique open
orbit in $\R\mathbb P^{n-1}$. It follows that $U\!\cdot G = G$ has a unique open orbit $\O$ in $Y=X_n\times
X_n\times \R\mathbb P^{n-1}$, so (\ref{hypo1}c) is also satisfied.

For the second example, $G$ is again $SL^{\pm 1}(n,\R)$ with the standard $\Q$-structure, $H = SL^{\pm 1}(2n,\R)
\times G$, also with the standard $\Q$-structure,
\begin{equation}
\label{ex2} G\, \ni\, g\ \ \mapsto \ \ \left( \begin{pmatrix}  g & 0
\\ 0 & g \end{pmatrix} \ , \ g \ \right)\ \in \ H \ = \ SL^{\pm 1}(2n,\R)\times G
\end{equation}
describes the embedding of $G$ into $H$, and $Y = X_{2n}\times \R\mathbb P^{n-1}$, which is a generalized flag
variety for $H$. The abelian subgroup
\begin{equation}
\label{ex2a} U \ \ = \ \ \left\{ \left. \left( \ \begin{pmatrix}  e
& c \\ 0 & e \end{pmatrix}\ , \ e \ \right) \ \right| \ c \in
M(n\times n,\R)\ \right\}
\end{equation}
of $H$ is unipotent, defined over $\Q$, and normalized by $G$. It can also be described as the unipotent radical of
the $(n,n)$ parabolic in $SL^{\pm 1}(2n,\R)$, hence as the unipotent radical of a $\,\Q$-parabolic subgroup of $H$.
Again the semisimplicity of $G$ implies b). The conjugation action of $G=SL^{\pm 1}(n,\R)$ on $M(n\times n,\R)$
preserves the Euclidean measure, so d) is satisfied. The assertion c) can be reduced to the corresponding assertion
about the first example: when we embed $X_n\times X_n$ diagonally into $X_{2n}$, and correspondingly $X_n\times
X_n\times \R\mathbb P^{n-1}$ into $Y$, the $U$-translates of $X_n\times X_n\times \R\mathbb P^{n-1}$ sweep out a
dense open subset of $Y$; thus, since $G$ has an open orbit in $X_n\times X_n\times \R\mathbb P^{n-1}$, $\,U\!\cdot
G = G\!\cdot U\,$ does have an open orbit $\O \subset Y$.

We need to briefly recall the notion of an automorphic distribution. Let $\G_H \subset H$ denote an arithmetic
subgroup. Then $H$ acts unitarily on $L^2(\G_H\backslash H)$, via right translation. An automorphic representation
consists of an irreducible unitary representation $(\pi,V)$ of $H$, together with a $H$-invariant embedding
\begin{equation}
\label{autorep} j\, :\, V \ \hookrightarrow \ L^2(\G_H\backslash H)\,.
\end{equation}
If $v\in V$ is a $C^\infty$ vector\begin{footnote}{i.e., a vector such that $H\ni h \mapsto \pi(h)v$ is a
$C^\infty$ $V$-valued function on $H$.}\end{footnote}, then $j(v)$ is a $\G_H$-invariant $C^\infty$ function on
$H$. That makes evaluation of $j(v)$ at $e\in G$ meaningful. The linear map
\begin{equation}
\label{autodist1}
\tau_j\,:\, V^\infty\, \longrightarrow\, \C\,,\ \ \ \tau_j(v)\,=\,j(v)(e)\,,
\end{equation}
defined on the space of $C^\infty$ vectors $V^\infty$,\,\ is continuous with respect to the intrinsic topology on
$V^{\infty}$,\,\ and is also $\G_H$-invariant. Thus $\tau_j$ can be regarded as a $\G_H$-invariant distribution
vector for the dual representation $(\pi',V')$,
\begin{equation}
\label{autodist2}
\tau_j \,\in\, ((V')^{-\infty})^{\G_H}.
\end{equation}
We refer to $\tau_j$ as the automorphic distribution associated to the automorphic representation (\ref{autorep}).
It completely determines $j$, since $V^\infty$ is dense in $V$ and $j(v)(h)=(\pi(h)j(v))(e)=\langle
\tau_j,\pi(h)v\rangle$ for all $v\in V^\infty$, $\,h\in H$.

To simplify the notation, we drop the subscript $j$, and we interchange the roles of $(\pi,V)$ and its dual
$(\pi',V')$; from now on,
\begin{equation}
\label{autodist3}
\tau \,\in\, (V^{-\infty})^{\G_H}.
\end{equation}
We shall also relax the assumption that $\tau$ corresponds to an irreducible subspace of $L^2(\G_H\backslash H)$,
as in (\rangeref{autorep}{autodist1}): in the discussion that follows, $\tau$ will denote an arbitrary
$\G_H$-invariant distribution vector for an admissible representation $(\pi,V)$ of finite length, on a reflexive
Banach space\begin{footnote}{For an expository discussion of these matters see \cite[\S
5 and appendix]{rapiddecay}.}\end{footnote}. In particular, $\tau$ may denote the automorphic distribution arising from an
Eisenstein series.

We return to the situation of a generalized real flag variety $Y\cong H/P$, with $P\subset H$ parabolic, but not
necessarily defined over $\,\Q\,$.\,\ To any finite dimensional complex representation
\begin{equation}
\label{erep}
\mu \,: P \ \rightarrow \ GL(E)
\end{equation}
we associate the $H$-invariant vector bundle $\mathcal E \to Y$ modeled on $E$. In other words, $\mathcal E$ is a
vector bundle to which the action of $H$ on $Y$ lifts, and whose fiber at the identity coset $y_0 \cong eP$ is
isomorphic to $E$ as a $P$-module. Tracing back through the definition, one obtains the description\vspace{-4pt}
\begin{equation}
\label{evb1}
C^\infty(Y,\E) \cong \{\,f : H \overset{C^\infty}{\longrightarrow} E \,\mid \, f(hp) = \mu(p^{-1})f(h)\,\ \text{for $h\!\in \!H,\, p\! \in \! P$}\,\}\vspace{2pt}
\end{equation}
of its space of $C^\infty$ sections. This isomorphism relates the action of $H$ on the space of smooth sections
$C^\infty(Y,\E)$, by means of the structure of equivariant vector bundle, to the action, via left translation, on
the space of smooth $E$-valued functions on $H$.

We follow the convention of using the term ``distribution" in the sense of ``generalized function": scalar valued
distributions are dual to smooth measures, and thus continuous functions can be viewed as distributions. In
complete analogy to (\ref{evb1}), there exists a natural, $H$-invariant isomorphism
\begin{equation}
\label{evb2}
C^{-\infty}(Y,\E) \cong \{\sigma : H \overset{C^{-\infty}}{\longrightarrow} E \,\mid \, \sigma(hp) = \mu(p^{-1})\sigma(h)\,\, \text{for $h\!\in \!H,\ p\! \in \! P$}\}
\end{equation}
between the space $C^{-\infty}(Y,\E)$ of sections of $\,\E$ with distribution coefficients and the space of
$E$-valued distributions on $H$ which transform in the specified manner under right translation by elements of $H$.
The identity $\sigma(hp) = \mu(p^{-1})\sigma(h)$ has only symbolic meaning, of course, since distributions cannot
be evaluated at points.

When one puts a Riemannian metric on $Y$ and a hermitian metric on the vector bundle $\,\E$, it makes sense to
consider the space of all $L^2$ sections of $\,\E \to Y$. That is a Hilbert space, on which $H$ acts continuously,
but generally not unitarily. The resulting representation of $H$ is known to be admissible, of finite length. Its
spaces of $C^\infty$ and distribution vectors are naturally isomorphic to the spaces (\ref{evb1}) and (\ref{evb2}),
respectively. From now on we keep fixed an arithmetic subgroup $\G_H \subset H$ and a particular
\begin{equation}
\label{hypo6}
\tau\,\in\, C^{-\infty}(Y,\E)^{\G_H}\,.
\end{equation}
Then $\tau$ is an automorphic distribution in the sense of our earlier discussion. According to the Casselman
embedding theorem \cite{Casselman:1980} and results of Casselman-Wallach \cite[\S 11]{wallach}, the space of
distribution vectors for any irreducible unitary representation can be embedded into the space of distribution
vectors (\ref{evb2}), with appropriate choices of $Y$ and $\,\E$. In that sense, all automorphic distributions
$\tau_j$ that encode irreducible subspaces of $L^2(\G_H\backslash H)$ as in (\rangeref{autorep}{autodist1}) can be
realized as in (\ref{hypo6}); for details see \cites{voronoi,rankin-selberg}.\footnote{In those references, we
consider only line bundles, rather than vector bundles $\mathcal E\to Y$. For linear groups, line bundles suffice.
In this paper we also consider nonlinear groups. The Cartan subgroups of a non-linear group $H$ need not be
abelian, which necessitates working with representations parabolically induced from finite dimensional
representations of dimension greater than one~-- or in geometric terms, with vector bundles rather than line
bundles.}

We denote the action of $H$ on sections of $\,\E\to Y$, and also on functions and distributions on $H$, by the
letter $\ell$, for left translation. In analogy to the case of automorphic forms, the automorphic distribution
$\tau$ is said to be cuspidal if
\begin{equation}
\label{cuspidal1} \int_{N/(\G_H\cap N)}\ell(n)\tau\,dn\ = \ 0
\end{equation}
for any unipotent subgroup $\{e\} \neq N\subset H$ which arises as the unipotent radical of a $\Q$-parabolic
subgroup $P\subset H$. We also need a weaker notion. Suppose
\begin{equation}
\label{cuspidal1.1} H=H_1\times H_2\ \ \ \text{(in the category of groups defined over $\Q$)}
\end{equation}
can be expressed as the product of two reductive, non-abelian factors. We shall then say that $\tau$ is cuspidal
with respect to the factor $H_1$ if the integral (\ref{cuspidal1}) vanishes for any unipotent subgroup $\{e\} \neq
N\subset H_1$ which arises as the unipotent radical of a $\Q$-parabolic subgroup $P\subset H_1$. That will be a
potential hypothesis in the statement of our main theorem, but under the additional assumption that
\begin{equation}
\label{cuspidal1.11}
\text{the projection of $G\subset H$ into $H_j$ has a finite kernel,}
\end{equation}
for both $j=1,\, 2$.

Recall the conditions (\ref{hypo1}). We choose a base point $\o \in \O$ and let $(UG)_\o$ denote the isotropy
subgroup of $\,U\!\cdot G\,$ at $\o$\,. The surjective map
\begin{equation}
\label{hypo1.1} p\,:\, U\!\cdot G \ \longrightarrow \ \O\,, \ \ \ \ p(u\,g)\,=\,u\,g\!\cdot \! \o\,,
\end{equation}
induces a $(U\!\cdot G)$-equivariant identification $\O \simeq (U\!\cdot G)/(UG)_\o$. Distribution sections of the
$H$-equivariant vector bundle $\,\E \to Y$ can be restricted to the open subset $\O\subset Y$, and then pulled back
from $\,\O \simeq (U\!\cdot G)/(UG)_\o\,$ to $\,U\!\cdot G\,$ via $p$. We choose $s\in H$ so that
\begin{equation}
\label{sdef}
\o \, = \, sP  \in H/P\, \cong \, Y, \, \ \text{and consequently}\, \ (UG)_\o \, = \, (U\!\cdot G) \cap sPs^{-1}.
\end{equation}
The $(U\!\cdot G)$-equivariant vector bundle $\,\E|_\O \to \O \cong (U\!\cdot G)/(UG)_\o$ is attached to the
representation
\begin{equation}
\label{erep1}
(U G)_\o \,\ni \, v \,\ \mapsto \,\ \mu(s^{-1}vs) \,\in\, GL(E)\,,
\end{equation}
hence, in analogy to (\ref{evb2}),
\begin{equation}
\label{hypo1.4}
C^{-\infty}(\O,\E) \, \simeq \, \{\,\sigma : (U\!\cdot G) \overset{C^{-\infty}}{\longrightarrow} E\,\mid \, r(v)\sigma = \mu(s^{-1}v^{-1}s)\sigma\ \,\text{for}\,\ v\in (UG)_\o \, \}\,.
\end{equation}
Here $r$ denotes the right translation action on possibly vector valued functions and its natural extension to
distributions.

Even if the representation (\ref{erep}) is irreducible, its restriction to $s^{-1}(U G)_\o s$ may well have a
trivial quotient. We suppose that is the case\footnote{This hypothesis is necessary to construct the distributional
pairing. Geometrically it means that the restriction of $\mathcal E$ to the open orbit $\mathcal O$ has a $U\!\cdot
G$-equivariantly trivial rank one quotient bundle. When we are dealing with a line bundle $\mathcal L$ rather than
a vector bundle $\mathcal E$~-- cf. the previous footnote~-- it reduces to the assumption that the restriction of
$\mathcal L$ to $\mathcal O$ is $U\!\cdot G$-equivariantly trivial. In some applications, such as \cite{extsquare},
it is vacuously satisfied.}, and fix
\begin{equation}
\label{erep2}
q\, : \, E \ \longrightarrow \ \C \ \ \ \ \ \text{\big(\,$(s^{-1}(U G)_\o s)$-equivariant projection\big)}.
\end{equation}
Composing, from right to left, restriction of sections of $\,\E$ from $Y$ to $\O$, pullback from $\O$ to $U\!\cdot
G$, and the projection $q$ from $E$ to $\C$, we obtain a map $\widetilde p^{\,*}$ from $C^{-\infty}(Y,\E)$ to
scalar valued distributions on $U\!\cdot G$,
\begin{equation}
\label{hypo1.6}
\widetilde p^{\,*}\, : \,  C^{-\infty}(Y,\E) \ \longrightarrow C^{-\infty}(U\!\cdot G)\,.
\end{equation}
In terms of the isomorphisms (\ref{evb2}) and (\ref{hypo1.4}), is given by the explicit formula
\begin{equation}
\label{hypo1.66}
\big(\widetilde p^{\,*}\sigma\big)(ug)\ = \ q\big(\sigma(ugs)\big)\,.
\end{equation}
Then $\widetilde p^{\,*}$ is $(U\!\cdot G)$-equivariant by construction. In particular, it maps $\G_H$-invariant
distribution sections, such as $\tau$, to scalar distributions invariant under both $\G_H\cap U$ and $\G_H \cap G$.

To continue with the hypotheses of our main theorem, we note that the two subgroups
\begin{equation}
\label{hypo2}
\G_G \ = \ \G_H \cap G \ \subset \ G\ \ \ \ \text{and} \ \ \ \ \G_U \ = \ \G_H \cap U \ \subset \ U
\end{equation}
are arithmetic. Since $G$ normalizes $U$, $\G_G$ normalizes $\G_U$; moreover
\begin{equation}
\label{hypo3}
\G_U \backslash U\, \ \text{is compact},
\end{equation}
as is the quotient of any unipotent linear group over $\Q$ modulo an arithmetic subgroup. For any character
\begin{equation}
\begin{aligned}
\label{hypo4}
&\psi \, : \, U \ \rightarrow \ \C^* \ \ \ \text{such that}
\\
&\ \ \ \ \ \text{a)}\ \, \psi \equiv 1\ \ \text{on}\ \ \G_U\,,\
\ \text{and}
\\
&\ \ \ \ \ \text{b)}\ \, \psi(gug^{-1}) = \psi(u)\ \ \text{for all}\,\ g\in G\,,\,\ u \in U\,,
\end{aligned}
\end{equation}
we define the averaging operator
\begin{equation}
\begin{aligned}
\label{hypo5} &A_\psi\, :\, C^{-\infty}(Y,\E)^{\G_H}\
\longrightarrow\ C^{-\infty}(Y,\E)^{\G_U\cdot \G_G}\,,
\\
&\qquad\qquad A_\psi (\sigma)\ \ =\ \frac 1 {\vol(U/\G_U)}\int_{U/\G_U} \psi(u)\,\ell(u)\sigma\,du\,.
\end{aligned}
\end{equation}
This makes sense because of (\ref{hypo1}d) and (\rangeref{hypo3}{hypo4}). We shall briefly comment on the analytic
aspects of the definition in section \ref{sec:smoothing}.

\begin{thm}\label{mainthm}
Under the hypotheses just stated, for any $\phi \in C_c^\infty(G)$,
\[
F_{\tau,q,\psi,\phi}\, =_{\text{def}}\, \int_G
\phi(g)\,r(g)\,\widetilde p^{\,*}\!\left(A_\psi \tau\right) dg
\]
is a left $\G_G$-invariant smooth function on $\,U\!\cdot G$, translating under the right action of $U$ on
$U\!\cdot G$ according to the character $\psi$, and under the left action of $U$ according to the complex conjugate
character $\overline{\psi}$. In particular the restriction of $F_{\tau,q,\psi,\phi}$ to $G$ lies in
$C^\infty(\G_G\backslash G)$. If $\tau \in C^{-\infty}(Y,\E)^{\G_H}$ is cuspidal or is at least cuspidal with
respect to one factor $H_j$ in the setting (\rangeref{cuspidal1.1}{cuspidal1.11}), this restricted function decays
rapidly along all the cusps of $\G_G\backslash G$, and the integral
\[
\int_{\G_G\backslash G} F_{\tau,q,\psi,\phi}(g)\,dg
\]
converges absolutely. The value of the integral does not depend on the choice of $\phi$, provided $\phi$ satisfies
the normalizing condition $\int_G \phi(g)\,dg=1$. If $\tau$ depends holomorphically on a
parameter\footnote{Holomorphic dependence is to be taken in the weak sense: a distribution $\sigma_s$ depends
holomorphically on the parameter $s$ if the integral of $\sigma_s$ against any test function depends
holomorphically on $s$. Eisenstein distributions~-- the analogues of Eisenstein series in the context of
automorphic distributions~-- are typical examples of distributions depending holomorphically on a parameter.}, then
the value of the integral also depends holomorphically on that parameter.
\end{thm}

Let us re-state the theorem informally, in more suggestive terms.  First of all, the action of $U\cdot G$ on its
open orbit $\mathcal O$ allows us to think of $\widetilde p^{\,*}\!\left(A_\psi \tau\right)\in C^{-\infty}(U\!\cdot
G)$ as
\begin{equation}\label{newexpl1}
  \widetilde p^{\,*}(A_\psi \tau)(ug) \ \ = \ \ \frac 1 {\vol(\G_U\backslash U)} \, \int_{\G_U\backslash U} \psi(u_1)\i\,q(\tau(u_1ugs))\,du_1\,.
\end{equation}
This distribution transforms by the character $\psi$ in the $u$-variable, and thus has smooth dependence on $u\in
U$.  It therefore makes sense to speak of its restriction to the $G$ factor, in which it is formally automorphic
under $\G_G$.  Suppose momentarily that $\G_G$ is cocompact in $G$.  Then the integral of this restriction  over
$\G_G\backslash G$,
\begin{equation}\label{newexpl2}
    \frac 1 {\vol(\G_U\backslash U)} \, \int_{\G_G\backslash G} \int_{\G_U\backslash U}
    \psi(u)\i\,q(\tau(ugs))\,du\,dg\,,
\end{equation}
makes sense as the integral of a distribution against the  smooth measure $dg$ over a compact manifold.  Right
translating  $g$ by $g_1\in G$ does not change the value of (\ref{newexpl2}), nor does integrating the resulting
integral over $g_1$ against a smooth function $\phi$ of compact support and total integral 1.  In other words, the
integral of the function
\begin{equation}\label{newexpl3}
    F_{\tau,q,\psi,\phi}(g) \ \ = \ \  \frac 1 {\vol(\G_U\backslash U)} \,\int_{G}  \int_{\G_U\backslash U}
    \psi(u)\i\,q(\tau(ugg_1s))\,\phi(g_1)\,du\,dg
\end{equation}
over $\G_G\backslash G$ --  the integral in the statement of the theorem -- equals (\ref{newexpl2}), again assuming
that $\G_G$ is cocompact in $G$.  In the noncompact setting (where nearly all our applications lie), we cannot
directly make sense of (\ref{newexpl2}), but rather use (\ref{newexpl3}) to define pairing in the statement of the
theorem.

In applications of  theorem~\ref{mainthm}, $H$ typically factors as a product $H = H_1 \! \times \! H_2\,$, and
correspondingly $\,\tau = \tau_1 \!\cdot \!\tau_2\,$  as the product of two automorphic distributions, one for each
of the factors $H_j$. One can then think of the integral of $ F_{\tau,q,\psi,\phi}$ in the theorem as defining a
pairing between the two automorphic distributions. That is the reason for the title of our paper. Often one of the
factors $\tau_j$ is an Eisenstein series -- not a classical Eisenstein series, of course, but its distribution
version. Nonetheless one can ``unfold", just as one does in the traditional Rankin-Selberg approach. This can be
carried out in one of several ways, depending on the particular application -- see, for example,
\cites{rankin-selberg,extsquare}. For that reason, we shall treat the unfolding not here, but in papers in which we
apply theorem \ref{mainthm}.

We conclude with some explicit examples of pairings that are explained in further detail in \cite{rankin-selberg}.
Recall that $s\in H$ was introduced in (\ref{sdef}) in order to identify the quotient $(U\!\cdot G)/(UG)_\o\,$ with
its open orbit $\mathcal O$.  It is a straightforward matter to make this explicit in any particular example. For
instance, let us consider the open orbit of $G=SL^{\pm 1}(2,\R)$ on $\RP^{1}\times \RP^1\times \RP^1$ that was
considered above in the earlier examples.  If $f_1=\ttwo1001$, $f_2=\ttwo1101$, and $f_3=\ttwo{0}{-1}{1}{0}$, the
diagonal action of $G$ on $H=G\times G\times G$ gives an open orbit  in the flag variety $\RP^{1}\times \RP^1\times
\RP^1$ with basepoint represented by $s=(f_1,f_2,f_3)\in H$.  Thus if $\tau_1$, $\tau_2$, and $\tau_3$ are
automorphic distributions for $G$, the integral defined in theorem~\ref{mainthm} is equal to
\begin{equation}\label{newexpl4}
    \int_{\G_G\backslash G}\int_G \tau_1(g h f_1)\,\tau_2(ghf_2)\,\tau_3(ghf_3)\,\phi(h)\,dh\,dg\,,
\end{equation}
and is related to the Rankin-Selberg $L$-function of $\tau_1 \otimes \tau_2$ when $\tau_3$ is an Eisenstein series
distribution. As an example with a nontrivial unipotent integration, let us return to (\ref{ex2})-(\ref{ex2a}) and
let $\tau_1$ be an automorphic distribution for $SL^{\pm 1}(4,\R)$, and $\tau_2$ be an automorphic distribution for
$G=SL^{\pm 1}(2,\R)$.  The integral defined in theorem~\ref{mainthm} in this case is
\begin{equation}\label{newexpl5}
\frac 1 {\vol(\G_U\backslash U)} \, \int_{\G_G\backslash G} \int_{G}  \int_{\G_U\backslash U}
\psi(u)\i\,\tau_1\(u\ttwo{ghf_1}{}{}{ghf_2}\)\,\tau_2(ghf_3)\,du\,dh\,dg\,,
\end{equation}
and is related to the exterior square $L$-function of $\tau_1$ when $\tau_2$ is an Eisenstein series distribution.

\section{Smoothing the integrand}\label{sec:smoothing}

In this section we will first prove the invariance and smoothness properties asserted in theorem~\ref{mainthm}, and
then give a construction of $F_{\tau,q,\psi,\phi}$ as an integral of an automorphic form on $H$ in
lemma~\ref{lem:smoothing}.  The integral is then used in the following section to prove the remaining properties
asserted in the theorem.

We continue with the hypotheses and notation of the previous section, and return briefly to the definition
(\ref{hypo5}) of the averaging operator $A_\psi\,$. Like the space of distributions on any $C^\infty$ manifold,
$C^{-\infty}(Y,\E)$ has an intrinsic topology as complete, locally convex, Hausdorff topological vector space.
These are precisely the properties one needs to define the integral of a conti\-nuous, compactly supported function
with values in a topological vector space. The action of $H$ on $C^{-\infty}(Y,\E)$ is conti\-nuous by general
principles. It follows that $u \mapsto \psi(u) \ell(u)\tau$, for $\tau\in C^{-\infty}(Y,\E)^{\G_H}$, can be
regarded as a conti\-nuous function, defined on the compact manifold $U/\G_U$, with values in the topological
vector space $C^{-\infty}(Y,\E)$. As such it has a well defined integral, which represents $A_\psi\tau$.

Let us first examine the smoothness and the invariance properties asserted by the theorem. For $u\in U$,
\begin{equation}
\label{smoothing1}
\begin{aligned}
\ell(u)(A_\psi\tau)\ &=\ \frac 1 {\vol(U/\G_U)}\int_{U/\G_U}\psi(v)\ell(uv)\tau\,dv
\\
& = \ \frac 1 {\vol(U/\G_U)}\int_{U/\G_U}\psi(u^{-1}v)\ell(v)\tau\,dv \ = \ \overline{\psi(u)}\,A_\psi\tau\,,
\end{aligned}
\end{equation}
and therefore, in view of the $(U\!\cdot G)$-equivariance of $\,\widetilde p^{\,*}$,
\begin{equation}
\label{smoothing2}
\ell(u)\,\widetilde p^{\,*}(A_\psi\tau)\ = \ \overline{\psi(u)}\,\,\widetilde p^{\,*}( A_\psi\tau)\,.
\end{equation}
In terms of the identification $U\!\cdot  G \simeq U\times G$, the distribution $\,\widetilde p^{\,*} (
A_\psi\tau)$ transforms according to the character ${\psi}$ in the $U$-variable. Without loss of information, we
can set the $U$-variable equal to the identity and -- temporarily~-- regard $\,\widetilde p^{\,*} ( A_\psi\tau)$ as
distribution on $G$. Convolution with a compactly supported smooth function turns any distribution on $G$ into a
smooth function. We then put the $U$-variable back in and conclude:
\begin{equation}
\label{smoothing3}
F_{\tau,q,\psi,\phi}\ =\ \int_G \phi(g)\,r(g)\, \widetilde p^{\,*}(A_\psi\tau)\, dg\ \ \ \text{is a $C^\infty$ function on $\,U\!\cdot G$\,},
\end{equation}
for every $\phi\in C^\infty(G)$. Convolution on the right commutes with left translation. Hence, for $u\in U$,
(\ref{smoothing2}) implies
\begin{equation}
\label{smoothing5} \ell(u)F_{\tau,q,\psi,\phi}\ =\
\overline{\psi(u)}F_{\tau,q,\psi,\phi}\,.
\end{equation}
Now suppose $g\in G$, $u, u_1\in U$. Then, in view of (\ref{hypo4}b) and (\ref{smoothing5}),
\begin{equation}
\label{smoothing6}
\begin{aligned}
\!\! (r(u)F_{\tau,q,\psi,\phi})(u_1g) &= F_{\tau,q,\psi,\phi}(u_1gug^{-1}g) = (\ell((u_1gug^{-1})^{-1})F_{\tau,q,\psi,\phi})(g)
\\
&= \psi(u_1gug^{-1})F_{\tau,q,\psi,\phi}(g) =
\psi(u)\psi(u_1)F_{\tau,q,\psi,\phi}(g)
\\
&=
\psi(u)\left(\ell(u_1^{-1})F_{\tau,q,\psi,\phi}\right)(g) =
\psi(u)F_{\tau,q,\psi,\phi}(u_1g)\,,
\end{aligned}
\end{equation}
so $\,r(u)F_{\tau,q,\psi,\phi}=\psi(u)F_{\tau,q,\psi,\phi}\,$. Since $A_\psi \tau$ is $\G_G$-invariant on the left,
since $\,\widetilde p^{\,*}$ is $G$-equivariant, and since right and left translation commute,
\begin{equation}
\label{smoothing6.5}
\ell(\g)F_{\tau,q,\psi,\phi}\ =\ F_{\tau,q,\psi,\phi}\,,\ \ \text{for every $\g\in \G_G$}\,.
\end{equation}
These are the smoothness and the invariance properties of $F_{\tau,q,\psi,\phi}\,$, as asserted.

For the proof of rapid decay in the next section, we need to work with an auxiliary function
$\Phi_{\tau,q,\psi,\phi,\phi_U}$, whose definition involves the choice of some $\phi_U \in C^\infty_c(U)$, in
addition to $\tau$, $\psi$, and $\phi$. We impose the normalization conditions
\begin{equation}
\label{smoothing7}
\int_G \phi(g)\,dg \ = \ 1\,,\qquad \int_U \phi_U(u)\,du \ = \ 1\,,
\end{equation}
and regard $\,q\!\circ\! \tau\,$ as a scalar valued distribution on $H$, via (\ref{evb2}) and (\ref{erep2}). That
allows us to consider
\begin{equation}
\label{smoothing8}
\Phi_{\tau,q,\psi,\phi,\phi_U}\ = \ \int_G\int_U \phi(g)\, \phi_U(u)\,\overline{\psi(u)}\,r(u)\,r(gs)\,q\!\circ\!\tau\, du\,dg \,,
\end{equation}
as a scalar distribution on $H$. Recall that $s\in H$, as defined in (\ref{sdef}), relates the base points $y_0\in
Y$ and $\o\in\O$.

\begin{lem}
\label{lem:smoothing} $\Phi_{\tau,q,\psi,\phi,\phi_U}$ is a $\G_H$-automorphic form\begin{footnote}{a smooth
automorphic form, not required to transform finitely under the right action of a maximal compact subgroup, as is
often assumed.}\end{footnote} on $H$ -- i.e., a $\G_H$-in\-va\-riant $C^\infty$ function of uniformly moderate
growth, which transforms finitely under the action of the algebra of bi-invariant differential operators on $H$.
Moreover, for all $g\in G$,
\[
\frac 1 {\vol(U/\G_U)}\int_{U/\G_U} \psi(u)\left(\ell(u)\,\Phi_{\tau,q,\psi,\phi,\phi_U}\right)(g)\,du\ = \ F_{\tau,q,\psi,\phi}(g)\,.
\]
\end{lem}

\noindent {\bf Proof:} Let $\,\E^*\to Y$ denote the $H$-invariant vector bundle dual to $\,\E$, and
$\wedge^{\text{top}}T^*Y$ the top exterior power of the cotangent bundle of $Y$ -- i.e., the line bundle whose
$C^\infty$ sections are smooth measures on $Y$. We shall produce a smooth section $\omega$ of the tensor product,
\begin{equation}
\label{smoothing8.1}
\omega \in C^\infty(Y,\E^*\otimes \wedge^{\text{top}}T^*Y)\,,\ \ \text{such that}\ \ \Phi_{\tau,q,\psi,\phi,\phi_U}(h)\, = \int_Y \langle\,\ell(h^{-1})\tau \,,\, \omega\,\rangle \ ;
\end{equation}
here $\langle\,\ell(h^{-1})\tau \,,\, \omega\,\rangle$ denotes the contraction between $\omega \in
C^\infty(Y,\E^*\otimes \wedge^{\text{top}}T^*Y)$ and $\ell(h\i)\tau \in C^{-\infty}(Y,\E)$, resulting in a scalar
valued, distribution coefficient form of top degree on the compact manifold $Y$. As such, it can be integrated over
$Y$ against any smooth function, in particular the constant function $1$.

We assume the existence of $\omega$ for the moment. In section \ref{sec:mainthm} we mentioned that
$C^{-\infty}(Y,\E)$ is the space of distribution vectors $V^\infty$ for an admissible, representation $(\pi,V)$ of
$H$, of finite length, on a Hilbert space. Analogously its topological dual $C^\infty(Y,\E^*\otimes
\wedge^{\text{top}}T^*Y)$ can be regarded as the space of $C^\infty$ vectors $(V')^\infty$ for the dual
representation $(\pi',V')$. If we now let $\langle \ ,\ \rangle$ denote the pairing between $V^{-\infty}$ and
$(V')^\infty$, we can rewrite the equality in (\ref{smoothing8.1}) as
\begin{equation}
\label{smoothing8.2}
\Phi_{\tau,q,\psi,\phi,\phi_U}(h)\ = \ \langle\,\ell(h^{-1})\tau\,,\, \omega\,\rangle \,.
\end{equation}
Since $\,h\mapsto \ell(h^{-1})\tau\,$ is a smooth function on $H$, with values in the topological vector space
$\,V^{-\infty}$,\,\ and since the pairing is linear and continuous in each variable, the description
(\ref{smoothing8.2}) exhibits $\Phi_{\tau,q,\psi,\phi,\phi_U}$ as a $C^\infty$ function. Moreover, any function of
this type, with $\tau \in (V^{-\infty})^{\G_H}$, is a $C^\infty$ automorphic form \cite[(2.15)]{rapiddecay}.

We continue to assume the existence of $\omega$. We shall establish the identity stated at the end of the lemma.
Since it only involves the values of $\Phi_{\tau,q,\psi,\phi,\phi_U}$ on $U\!\cdot G$, we now regard this function
as defined on $U\cdot G$, rather than on $H$ as in (\ref{smoothing8}). Recall the definition of the averaging
operator $A_\psi$ in (\ref{hypo5}). Left and right translation commute, hence
\begin{equation}
\label{smoothing15}
\begin{aligned}
&\frac 1 {\vol(U/\G_U)}\int_{U/\G_U}
\psi(u)\left(\ell(u)\,\Phi_{\tau,q,\psi,\phi,\phi_U}\right)\,du\ =
\\
&\qquad = \ \int_U\int_G \phi(g)\,
\phi_U(u)\,\overline{\psi(u)}\,r(u)\,r(g)\,\widetilde p^{\,*}(A_\psi\tau)\, dg\,du
\\
&\qquad = \ \int_U
\phi_U(u)\,\overline{\psi(u)}\,r(u)F_{\tau,q,\psi,\phi}\,du
\\
&\qquad = \ \int_U \phi_U(u)\,|\psi(u)|^2\,F_{\tau,q,\psi,\phi}\,du \
= \ F_{\tau,q,\psi,\phi}\,\,.
\end{aligned}
\end{equation}
At the first step of this chain of equalities, we have used the definition (\ref{smoothing8}) of
$\Phi_{\tau,q,\psi,\phi,\phi_U}$, the definition (\ref{hypo5}) of $\,A_\psi\,$,\,\ as well as the characterization
of the map $\widetilde p^{\,*}$ in (\ref{hypo1.6}); at the second step we have used the definition of
$F_{\tau,q,\psi,\phi}$ in the statement of the main theorem; and finally, in the third step, the transformation law
(\ref{smoothing6}) and the normalization (\ref{smoothing7}) of $\phi_U$.

We shall construct $\omega$ first as a smooth, compactly supported, $\,\E$-valued measure on $\O$. Since $\O$ is
open in $Y$, we can then regard $\omega$ as an element of $C^\infty(Y,\E^*\otimes \wedge^{\text{top}}T^*Y)$. Recall
the definition, above (\ref{hypo1.1}), of $(UG)_\o$ as the isotropy subgroup of $\,U\!\cdot G\,$ at $\o =
s\!\cdot\!y_0 \cong sP\,$. Since $\,U\!\cdot G\,$ and $\,P\,$ are defined, respectively, over $\,\Q\,$ and over
$\,\R\,$,
\begin{equation}
\label{vdef}
(UG)_\o\ = \ (U\!\cdot G)\cap sPs^{-1}
\end{equation}
is an $\,\R$-subgroup of $\,H\,$. We define
\begin{equation}
\label{fdef}
f_0\in C^\infty_c(U\!\cdot G)\,,\ \ \ f_0(ug)\ = \ \phi(g)\,\phi_U(u)\,\overline{\psi(u)}\,.
\end{equation}
Note that $\,U\!\cdot G\,$ is unimodular because of (\ref{hypo1}d). Haar measure on this group can be described as
the product $\,du\,dg\,$ of the Haar measures on $U$ and $G$. The normalizations (\ref{smoothing7}) imply that
$\,f_0\,du\,dg\,$ is well-defined, independently of the scaling of $\,du\,$ and $\,dg\,$. We let $(\mu^*,E^*)$
denote the representation of $P$ dual to $(\mu,E)$, and choose $e^*\in E^*$ such that
\begin{equation}
\label{e*def}
q(e) \ = \ \langle \,e^*\,,\, e\,\rangle \ \ \ \ \text{for all}\ \ e \in E\,.
\end{equation}
Using this notation, we interpret (\ref{smoothing8}) as an identity between distributions on $H$,
\begin{equation}
\label{smoothing9.5}
\Phi_{\tau,q,\psi,\phi,\phi_U}(h) \ = \ \int_{U\cdot G}  \langle \,f_0(ug)\,e^*\,,\, \tau(hugs)\,\rangle \,du\,dg\,.
\end{equation}
The group $(UG)_\o$ may not be unimodular. We pick a left invariant Haar measure $dv\,$;\,\ the particular
normalization will turn out not to matter. Averaging first over $(UG)_\o$ and integrating the result over the
quotient $\,(U\!\cdot G)/(UG)_\o\,$, we find
\begin{equation}
\label{smoothing9.6}
\Phi_{\tau,q,\psi,\phi,\phi_U}(h) \, = \, \int_{(U\cdot G)/(UG)_\o}\int_{(UG)_\o} \!\! \langle \,f_0(ugv)\,e^*\,,\, \tau(hugvs)\,\rangle \, dv\ \frac{du\,dg}{dv}\ ;
\end{equation}
here $\,du\,dg/dv\,$ is regarded as smooth measure on $\,(U\!\cdot G)/(UG)_\o\,$ in the obvious manner.

All along we have used (\ref{evb2}) to regard $\tau$ as $E$-valued distribution on $H$, transforming under right
translation by $p\in P$ according to $\mu(p^{-1})$. Hence
\begin{equation}
\begin{aligned}
\label{smoothing9.7}
\langle \,f_0(ugv)\,e^*\,,\, \tau(hugvs)\,\rangle \, &= \, \langle \,f_0(ugv)\,e^*\,,\, \mu(s^{-1}v^{-1}s)\tau(hugs)\,\rangle
\\
&= \, \langle \,f_0(ugv)\, \mu^*(s^{-1}vs)e^*\,,\,\tau(hugs)\,\rangle  \,.
\end{aligned}
\end{equation}
The $C^\infty$ function $\,f : U\!\cdot G \to E^*\,$,\,\ defined by
\begin{equation}
\label{Fdef}
f(ug)\ = \ \int_{(UG)_0} f_0(ugv)\,\mu^*(s^{-1}vs)e^*\,dv\,,
\end{equation}
is well defined because the support of any left translate of $f_0$ intersects $\,(UG)_\o\,$ in a compact subset. As
direct consequence of the construction,
\begin{equation}
\label{Fdef2}
f(ugv)\ = \ \mu^*(s^{-1}v^{-1}s)f(ug)\,\ \ \ \ \text{for all $\,v\in (UG)_\o$}\,.
\end{equation}
That, in analogy to (\ref{hypo1.4}) with $C^\infty$ in place of $C^{-\infty}$, allows us to regard $f$ as a
$C^\infty$ section of $\,\E^* \to (U\!\cdot G)/(UG)_\o \cong \O\,$. The fact that $f_0$ has compact support in
$\,U\!\cdot G\,$ implies the compact support, modulo $(UG)_\o $, of $f$,
\begin{equation}
\label{Fdef3}
f \in C^\infty_c(\O,\E^*)\,.
\end{equation}
At this point we can reinterpret the outer integral in (\ref{smoothing9.6}) as the integral of the smooth,
compactly supported, $\,\E^*$-valued measure
\begin{equation}
\label{omegadef}
\omega \ \ = \ \ f \, \frac{du\,dg}{dv} \, \in \, C^\infty_c(\O,\E^*\otimes \wedge^{\text{top}}T^*Y)
\end{equation}
against the $\,\E$-valued distribution section $\,\ell(h^{-1})\tau \in C^\infty(Y,\E)$. In effect, we have verified
(\ref{smoothing8.1}), completing the proof of the lemma. \bx

To construct $\omega$ in the preceding proof, we averaged $f_0$ over $(UG)_0$ on the right, and $f_0$ involved the
function $\phi$ as a factor. Hence:

\begin{rem}
\label{rem:omega} If the support of $\phi$ is kept fixed, the dependence of $\omega$ on $\phi$ is bounded in $C^k$
norm, for every $k\in \N$.
\end{rem}

\section{Rapid decay}\label{sec:proofofmainthm}

In this section we finish the proof of theorem~\ref{mainthm}, using lemma~\ref{lem:smoothing} as a main tool. We
shall use the results of \cite{rapiddecay}, in particular theorem 4.7 of that paper. The notation there differs
from our current notation, since the two papers were written for different purposes. The role of the current
ambient group $H$, or of its factors $H_j$ in (\ref{cuspidal1.1}), is played by $G$ in \cite{rapiddecay}, and the
role of the unipotent group $U$ by $N_1$. The group $G$ in \cite{rapiddecay} satisfies exactly the same hypotheses
as our ambient group $H$. Since our current $G$ is reductive and defined over $\,\Q\,$,\,\ it can be extended to a
Levi component, defined over $\,\Q\,$,\,\ of the parabolic subgroup $P_U \subset H$. Thus we may assume that the
Langlands decomposition in (\ref{hypo1}b) is defined over $\Q$. That, in view of (\ref{hypo4}), makes $G$ a
reductive $\,\Q$-subgroup of the group denoted by $L$ in \cite[\S 4]{rapiddecay}. If $L_1 \subset L_2$ is an
inclusion of real reductive groups defined over $\,\Q\,$,\,\ any Siegel set in $L_1$ is contained in the
intersection of $L_1$ with an appropriately chosen Siegel set in $L_2$. Thus, for statements about rapid decay on
Siegel sets of automorphic forms on some ambient group, any decay statement for $L_2$ directly implies the
analogous statement for $L_1$. In our words, for our purposes it is really irrelevant whether the current $G$ is a
subgroup of the $L$ in \cite[\S 4]{rapiddecay}, or all of it.

Let us restate theorem 4.7 of \cite{rapiddecay} in our current notation. We consider the automorphic form
$\Phi_{\tau,q,\psi,\phi,\phi_U}$ on $H$ associated to a cuspidal automorphic distribution $\,\tau \in
C^{-\infty}(Y,\E)$ and $\omega \in C^\infty_c(\O,\E^*\otimes \wedge^{\text{top}}T^*Y)$ as in (\ref{smoothing8.1});
in the proof of lemma \ref{lem:smoothing}, $\omega$ was constructed in terms of $\,q,\psi,\phi,\phi_U\,$,\,\ of
course. Then for any Siegel set $\mathfrak S_G \subset G$, any $c$ in the inverse image of $G_\Q$ in $G$~-- cf. the
discussion at the beginning of section \ref{sec:mainthm} -- and any $n \in \N$, there exists a constant $C>0$ such
that
\begin{equation}
\label{decay1}
u \in U\,,\ \ g \in \mathfrak S_G \ \ \ \Longrightarrow \ \ \ |\Phi_{\tau,q,\psi,\phi,\phi_U}(u\,c\,g)| \, \leq \, C\|g\|^{-n}\,.
\end{equation}
When $F_{\tau,q,\psi,\phi}(u\,c\,g)$ is defined in terms of $\Phi_{\tau,q,\psi,\phi,\phi_U}$ by the formula in
lemma \ref{lem:smoothing}, the above estimate directly implies
\begin{equation}
\label{decay2}
g \in \mathfrak S_G \ \ \ \Longrightarrow \ \ \ |F_{\tau,q,\psi,\phi}(u\,c\,g)| \, \leq \, C\|g\|^{-n}\,.
\end{equation}
A finite number of translates $c\,\mathfrak S_G$ cover $\G_G\backslash G$, so (\ref{decay2}) implies the
integrability of $F_{\tau,q,\psi,\phi}(u\,c\,g)$ over $\G_G\backslash G$, at least when $\tau$ is cuspidal.

The preceding argument must be modified when $\tau$ is only cuspidal with respect to one of the factors in a
factorization (\ref{cuspidal1.1})~-- say with respect to $H_1$ for definiteness. In that situation \cite[theorem
5.13]{rapiddecay} asserts the simultaneous rapid decay in the first factor and moderate growth in the second for
$\,\Phi_{\tau,q,\psi,\phi,\phi_U}(h_1h_2)\,$,\,\ provided $h_1$ varies in a Siegel set $\mathfrak S_1\subset H_1$
and $h_2$ over $H_2$. Rapid decay trumps moderate growth when $\,\Phi_{\tau,q,\psi,\phi,\phi_U}\,$ is restricted to
$U\!\cdot G\,$; this depends on the hypothesis (\ref{cuspidal1.11}), of course. Thus both (\ref{decay1}) and
(\ref{decay2}) remain valid even in the partially cuspidal case.

Only two assertions remain to be verified to complete the proof of theorem \ref{mainthm}. First of all, we need to
show that the integral
\begin{equation}
\label{smoothing16}
\int_{\G_G\backslash G} F_{\tau,q,\psi,\phi}(g)\,dg \ = \ \int_{\G_G\backslash G}\, \int_G \phi(g_1)\,\widetilde p^{\,*}(A_\psi\tau)(gg_1)\,dg_1\,dg
\end{equation}
depends only on the total integral of the smoothing function $\phi$, but not otherwise on the particular choice of
$\phi$. Since the measure on $\G_G\backslash G$ is $G$-invariant on the right, we can right translate the argument
$g$ by any fixed $g_2\in G$ without changing the value of the integral (\ref{smoothing16}). Equivalently, we can
replace $\phi$ by $\ell(g_2)\phi$ without affecting the integral. On the infinitesimal level this means that the
integral vanishes whenever $\phi = \ell(Z)\tilde\phi$ is the left infinitesimal translate of some $\tilde\phi\in
C^\infty_c(G)$ by some $Z\in\fg\,$. Thus the integral (\ref{smoothing16}) depends only on the image of $\phi$ in
the quotient
\begin{equation}
\label{smoothing17}
C_c^\infty(G)/\ell(\fg)C^\infty_c(G)\,.
\end{equation}
Since $G$ is parallelizable, this quotient represents the top cohomology group of the de Rham complex with compact
support. As such, it is naturally dual to $H_0(G,\C)$, and the duality is implemented by integration over the
various connected components of $G$. If $G$ is connected, we can conclude that
\begin{equation}
\label{smoothing18}
\ell(\fg)C^\infty_c(G)\ = \ \{\,\phi \in C_c^\infty(G)\, \mid \, \textstyle\int_G \phi(g)\,dg = 0\,\}\,,
\end{equation}
which implies the conclusion we want. The case of several connected components reduces readily to the connected
case; what matters here is that replacing $\phi$ by $\ell(g_2)\phi$, with $g_2 \in G$, does not affect the value of
the integral.

The preceding argument involves interchanging the order of differentiation and integration. To make this
legitimate, one needs to know that the absolute convergence of the integral of $F_{\tau,q,\psi,\phi}$ over
$\G_G\backslash G$ can be estimated~-- when $\tau$, $\psi$, and the support of $\phi$ are kept fixed~-- in terms of
some $C^k$ norm on $\phi$. Recall the remark \ref{rem:omega}. Since $\,F_{\tau,\psi,\phi}\,$ was expressed in terms
of the automorphic form $\,\Phi_{\tau,\psi,\phi,\phi_U}\,$,\,\ which involves $\phi$ via $\omega$ in the identity
(\ref{smoothing8.2}), we need to quantify the decay of $\,\Phi_{\tau,\psi,\phi,\phi_U}\,$ in terms of some $C^k$
norm of $\omega$, again assuming $\tau$ is fixed. In the notation of \cite[theorem 5.13]{rapiddecay}, $\omega$
corresponds to the vector $v$. The estimate in that theorem depends linearly on the $C^k$ norm of $v$, which
provides the necessary justification.

As was remarked earlier, for us holomorphic dependence of a distribution on a complex parameter means dependence in
the ``weak sense" -- i.e., the integral of the distribution against any test function depends holomorphically on
the parameter in question. According to (\ref{smoothing8.1}), the values of the function
$\Phi_{\tau,q,\psi,\phi,\phi_U}$ can be interpreted as the result of paring the distribution section $\tau$ of the
vector bundle $\mathcal E$ over the compact manifold $Y$ against a smooth measure with values in the dual vector
bundle. Thus, if $\tau$ depends holomorphically on a complex parameter, then so does the function
$\Phi_{\tau,q,\psi,\phi,\phi_U}$, and in view of lemma \ref{lem:smoothing}, also the function
$F_{\tau,q,\psi,\phi}$. This completes the proof of theorem \ref{mainthm}.

\begin{tabular}{lcl}
Stephen D. Miller                    & & Wilfried Schmid \\
Department of Mathematics         & & Department of Mathematics \\
Hill Center-Busch Campus          & & Harvard University \\
Rutgers, The State University of New Jersey             & & Cambridge, MA 02138 \\
 110 Frelinghuysen Rd             & & {\tt schmid@math.harvard.edu}\\
 Piscataway, NJ 08854-8019\\
 {\tt miller@math.rutgers.edu}  \\

\end{tabular}

\end{document}